\documentclass[12pt]{article}

\usepackage{mathptm}    
\usepackage{latexsym}   
\usepackage{amssymb}    
\usepackage{amsmath}    
\usepackage{amsthm}     

\newtheorem{theorem}{Theorem}[section]
\newtheorem{lemma}[theorem]{Lemma}
\newtheorem{proposition}[theorem]{Proposition}
\newtheorem{corollary}[theorem]{Corollary}

\newenvironment{prf}[1]{\trivlist
\item[\hskip
\labelsep{\it #1.\hspace*{.3em}}]}{
\endtrivlist}

\newtheorem{predefinition}[theorem]{Definition}

\newtheorem{preremark}[theorem]{Remark}
\newenvironment{remark}{\begin{preremark}\rm}{\end{preremark}}
\newtheorem{prenotation}[theorem]{Notation}

\newtheorem{preexample}[theorem]{Example}

\newtheorem{preclaim}[theorem]{Claim}

\newtheorem{prequestion}[theorem]{Question}
\newenvironment{question}{\begin{prequestion}\rm}{\end{prequestion}}
 \makeatletter
\def\emppsubsection{\@startsection{subsection}{2}{\z@}{-3.25ex plus -1ex minus -.2ex}{-1em}}

\makeatother

\newcommand \CA {{\cal A}}

\newcommand \CM {{\cal M}}

\newcommand \CH {{\cal H}}

\newcommand \PP {{\mathbb P}^1}
\newcommand \Aa {{\mathbb A}^1}
\newcommand \ZZ {{\mathbb Z}}

\newcommand  \FF {{\mathbb F}}
\newcommand \CC {{\mathbb C}}

\newcommand \GG {{\mathbb G}}

\newcommand \Aut {\mathop{\rm Aut}}

\newcommand \Det {\mathop{\rm Det}}
\newcommand \dime {\mathop{\rm dim}}

\newcommand \gen {\mathop{\rm genus}}

\newcommand \Hom {\mathop{\rm Hom}}

\newcommand \Jac {\mathop{\rm Jac}}

\newcommand \Spec {\mathop{\rm Spec}}

\newcommand \crys{{{\rm crys}}}

\title{Hyperelliptic curves with prescribed $p$-torsion}
\author{Darren Glass \& Rachel Pries
\footnote{The authors were partially supported by NSF VIGRE grant
DMS-98-10750.}}
\date{}

\begin{document}
\maketitle

\begin{abstract} \noindent In this paper, we show that there exist families of
curves (defined over an algebraically closed field $k$ of characteristic
$p >2$) whose Jacobians have interesting $p$-torsion.
For example, for every $0 \leq f \leq g$,
we find the dimension of the locus of hyperelliptic curves of genus $g$ with $p$-rank at most $f$.
We also produce families of curves so that the $p$-torsion of the Jacobian of each fibre
contains multiple copies of the group scheme $\alpha_p$.
The method is to study curves which admit an action by $(\ZZ/2)^n$ so that the quotient is a projective line.
As a result, some of these families intersect the hyperelliptic locus $\CH_g$.
\end{abstract}

\section{Introduction}

When investigating abelian varieties defined over an algebraically
closed field $k$ of characteristic $p$, it is natural to study the
invariants related to their $p$-torsion such as their $p$-rank or
$a$-number.  Such invariants are well-understood and have been
used to define stratifications of the moduli space $\CA_g$ of
principally polarized abelian varieties of dimension $g$.  There
is a deep interest in understanding whether the Torelli locus
intersects such strata in $\CA_g$. More generally, one can ask for
the dimension of the intersection of these strata with the image
of the moduli spaces $\CM_g$ or $\CH_g$ under the Torelli map.
In this paper, we show that the Torelli locus intersects several of
these strata by producing families of curves so that the
$p$-torsion of the Jacobian of each fibre contains certain group schemes.

Recall that the group scheme $\mu_p=\mu_{p,k}$ is the kernel of Frobenius on $\GG_m$
and the group scheme $\alpha_p=\alpha_{p,k}$ is the kernel of Frobenius on $\GG_a$.
As schemes, $\mu_p \simeq \Spec(k[x]/(x-1)^p)$ and $\alpha_p \simeq \Spec(k[x]/x^p)$ over $k$.
If $\Jac(X)$ is the Jacobian of a $k$-curve $X$,
the {\it $p$-rank} of $X$ is $\dime_{\FF_p} \Hom(\mu_p, \Jac(X))$
and the {\it $a$-number} of $X$ is $\dime_k \Hom(\alpha_p, \Jac(X))$.

Let $V_{g,f}$ denote the sublocus of $\overline{\CM}_g$ consisting of curves of genus $g$ with
$p$-rank at most $f$. For every $g$ and every $0 \leq f \leq g$,
the locus $V_{g,f}$ has codimension $g-f$ in $\overline{\CM}_g$,
\cite{FVdG}. In Section \ref{Sprank}, we use results from
\cite{FVdG} to prove that there exist smooth hyperelliptic curves
of genus $g$ with every possible $p$-rank $f$.

\bigskip

\noindent {\bf Theorem \ref{Thyperprk}.} For all $g \geq 1$ and all $0 \leq f \leq g$,
the locus $V_{g,f} \cap \CH_g$ is non-empty of
dimension $g-1+f$. In particular, there exists a smooth
hyperelliptic curve of genus $g$ and $p$-rank $f$.

\bigskip

Let $T_{g,a}$ denote the sublocus of $\overline{\CM}_g$ consisting of curves of genus $g$ with $a$-number at least $a$.
In Section \ref{Sanumb}, we show that $T_{g,a}$ is non-empty under certain conditions on $g$ and $a$
by producing curves $X$ so that $\Jac(X)[p]$ contains multiple copies of $\alpha_p$.
Let $\CH_{g,n}$ be the sublocus of the moduli space $\CM_g$ consisting
of smooth curves of genus $g$ which admit an action by $(\ZZ/2)^n$
so that the quotient is a projective line.

\bigskip

\noindent {\bf Corollary \ref{Pa=n}.} Suppose $n \geq 2$ and $p \geq
2n+1$. Suppose $g$ is such that $\CH_{g,n}$ is non-empty of
dimension at least $n+1$. Then the intersection $\CH_{g,n} \cap
T_{g,n}$ has codimension at most $n$ in $\CH_{g,n}$. In
particular, there exists a smooth curve of genus $g$ with
$a$-number at least $n$.

\bigskip

The dimension of the family in Corollary \ref{Pa=n} is at least $(g+2^n-1)/2^{n-2}-3-n$.
The precise numerical conditions for $g$ can be found in Section \ref{Sanumb}.
The main interest in this result is not only that certain group schemes
occur in the $p$-torsion of the Jacobians, but also that the
dimension of the families is large in comparison with the
dimension of $\CH_{g,n}$.

For small values of $n$, we further show that these families of curves intersect the
hyperelliptic locus $\CH_g$, resulting in the following corollaries.

\bigskip

\noindent {\bf Corollary \ref{Tanumb2}.}
Suppose $g \geq 2$ and $p \geq 5$.   There
exists a $(g-2)$-dimensional family of smooth hyperelliptic curves of genus
$g$ whose fibres have $a$-number $2$ and $p$-rank $g-2$.

\bigskip

\noindent {\bf Corollary \ref{Canum3}.}
Suppose $g \geq 5$ is odd and $p \geq 7$.   There exists a
$(g-5)/2$-dimensional family of smooth hyperelliptic curves of genus $g$
whose fibres have $a$-number at least $3$.

\bigskip

In Section \ref{Sptor}, we consider the problem of constructing Jacobians whose $p$-torsion
contains group schemes other than $\mu_p$ or $\alpha_p$.
We prove that for all $g \geq 2$ there exists a
smooth hyperelliptic curve of genus $g$ whose $p$-torsion contains the group scheme corresponding
to a supersingular non-superspecial abelian surface.
We describe this group scheme and its covariant Dieudonn\'e module in Section \ref{Sptor}.
It has $a$-number $1$ and $p$-rank $0$.

Our method for these results is to analyze the curves in the locus $\CH_{g,n}$ in terms of fibre
products of hyperelliptic curves.  In Section \ref{S2}, we extend results of Kani and Rosen
\cite{KR} to compare the $p$-torsion of the Jacobian of a curve $X$ in
$\CH_{g,n}$ to the $p$-torsion of the Jacobians of its $\ZZ/2\ZZ$-quotients
{\it up to isomorphism}.  We then use Yui's description of the branch locus of
a non-ordinary hyperelliptic curve, \cite{Y}.   In some cases, this reduces the study of the $p$-torsion of
the Jacobian of $X$ to the study of the intersection of some subvarieties in
the configuration space of branch points.  We consider this in Section \ref{S3}.


Throughout, $k$ is an algebraically closed field of characteristic $p >2$.
We assume $g \geq 1$ to avoid trivial cases.
Without further comment, we will speak of a {\it fibre} of a relative curve
when we mean a geometric fibre.

This paper led us to pose some open questions on this topic in
\cite{GP:questions}.  We would like to thank E. Goren for
suggesting the topic of this paper, E. Kani for help with
Proposition \ref{Pcohom}, and J. Achter, I. Bouw, F. Oort, and the
referee for helpful comments.

\section{Curves with prescribed $p$-rank} \label{Sprank}

We begin by considering the $p$-rank of Jacobians of hyperelliptic
curves. Recall that the $p$-rank, $\dime_{\FF_p} \Hom(\mu_p,
\Jac(X))$, of a $k$-curve $X$ is an integer between $0$ and its
genus $g$. The curve $X$ is said to be {\it ordinary} if it has
$p$-rank equal to $g$. In other words, $X$ is ordinary if
$\Jac(X)[p] \cong (\ZZ/p \oplus \mu_p)^g$. Let $V_{g,f}$ denote
the sublocus of $\overline{\CM}_g$ consisting of curves of genus
$g$ with $p$-rank at most $f$.

Consider the moduli space $\CH_g$ of smooth hyperelliptic curves of genus $g$
and its closure $\overline{\CH}_g$ in $\overline{\CM}_g$.
It is known that $\CH_g$ is affine of dimension $2g-1$.
Both $\CH_g$ and $\overline{\CH}_g$ are smooth algebraic stacks over $\ZZ[1/2]$
(for example, see \cite[Proposition 1]{W:thesis}).
Since $k$ is an algebraically closed field,
this fact implies that if two subvarieties of $\overline{\CH}_g$ intersect then the codimension of
their intersection is at most the sum of their codimensions.

The boundary $\overline{\CH}_g -\CH_g$ consists of components
$\Delta_0$ and $\Delta_i$ for integers $1 \leq i \leq g/2$. The
generic point of $\Delta_i$ corresponds to the isomorphism class
of a singular curve with two irreducible components $X_i$ and
$X_{g-i}$ intersecting in a node which we denote $P_i$. Here $X_i$
(resp.\ $X_{g-i}$) is a hyperelliptic curve of genus $i$ (resp.\
$g-i$) and the point $P_i$ is fixed by the hyperelliptic
involution on $X_i$ (resp.\ $X_{g-i}$). The generic point of
$\Delta_0$ corresponds to the isomorphism class of an irreducible
hyperelliptic curve $X_0'$ with a node. The normalization $X_0$ of
$X_0'$ is a hyperelliptic curve of genus $g-1$, and the inverse
image of the node in $X_0'$ consists of two distinct points in $X_0$ which
are exchanged by the hyperelliptic involution. Note that
$\Jac(X_0')$ is a semi-abelian variety and the toric part of its
$p$-torsion contains a copy of the group scheme $\mu_p$.  So
$\Delta_0 \cap V_{g,0}$ is empty in $\overline{\CM}_g$.

We first show that each component of $V_{g,0} \cap \overline{\CH}_g$ has dimension $g-1$.

\begin{proposition}\label{Pprk0}
The locus $V_{g,0} \cap \overline{\CH}_g$ is pure of codimension $g$ in $\overline{\CH}_g$.
\end{proposition}

\begin{proof}
We work by induction on $g$.
The statement is true in the case $g=1$ since the locus of supersingular
elliptic curves has dimension $0$.  Assume that the
statement is true for all $g'<g$, and consider any component $C_0$
of the intersection $V_{g,0} \cap \overline{\CH_g}$. By the purity
argument of \cite[1.6]{OdJ}, the codimension of $C_0$ in
$\overline{\CH}_g$ is at most $g$. Furthermore, $C_0$ intersects
the boundary of $\overline{\CH}_g$ because $\CH_g$ is affine.
Since $C_0$ does not intersect $\Delta_0$, it must intersect
$\Delta_i$ for some $1 \leq i \leq g/2$. We fix one such
$\Delta_i$ and consider the dimension of the intersection.

A curve corresponding to a point in the intersection of $C_0$
and $\Delta_i$ is formed from two hyperelliptic curves $X_i$ and
$X_{g-i}$ which must both have $p$-rank $0$. Thus $X_i$
corresponds to a point of $V_{i,0} \cap \overline{\CH}_i$ and
likewise $X_{g-i}$ corresponds to a point of $V_{g-i,0} \cap
\overline{\CH}_{g-i}$. By the inductive hypothesis, there is at
most an $i-1$ (resp.\ $g-i-1$) dimensional family of choices for
$X_i$ (resp.\ $X_{g-i}$). Since $X_i$ and $X_{g-i}$ intersect in a
unique point $P_i$, this point must be fixed under the
hyperelliptic involutions of the two curves.  Thus there are
only finitely many choices for the point $P_i$.
It follows that $\dime(C_0 \cap \Delta_i) \leq
(i-1)+(g-i-1)+0 = g-2$ and the codimension of $C_0 \cap \Delta_i$
in $\overline{\CH}_g$ is at least $g+1$.

We can deduce that ${\rm codim}(C_0 \cap \Delta_i) \leq {\rm
codim}(C_0) +1$ in $\overline{\CH}_g$ from the fact that
$\Delta_i$ has codimension $1$ in $\overline{\CH}_g$. This implies
that the codimension of $C_0$ in $\overline{\CH}_g$ is exactly $g$
and therefore that $V_{g,0} \cap \overline{\CH}_g$ is pure of
codimension $g$ in $\overline{\CH}_g$.
\end{proof}

Next we show that each component of $V_{g,f} \cap
\overline{\CH}_g$ has dimension $g-1+f$ (for $g \geq 1)$.

\begin{proposition}\label{Phyperprk}
The locus $V_{g,f} \cap \overline{\CH}_g$ is pure of codimension $g-f$ in $\overline{\CH}_g$.
\end{proposition}

\begin{proof}  By Proposition \ref{Pprk0}, we can suppose $f \geq 1$.
Consider a component $C_0$ of $V_{g,f} \cap \overline{\CH}_g$. By
\cite[1.6]{OdJ}, $C_0$ has codimension at most $g-f$ in
$\overline{\CH}_g$ and thus dimension at least $g-1+f$. Because $p
>2$, a complete subvariety of $\overline{\CH}_g - \Delta_0$ has
dimension at most $g-1$, by \cite[Lemma 2.6]{FVdG}. So $C_0$
intersects $\Delta_0$.

A point of $C_0 \cap \Delta_0$ corresponds to a curve $X_0'$
self-intersecting in a node $P_0$. The normalization $X_0$ of
$X_0'$ is a hyperelliptic curve of genus $g-1$. Since the toric
part of $\Jac(X_0')[p]$ contains a copy of the group scheme
$\mu_p$, this implies that the $p$-rank of $X_0$ is at most $f-1$.
So $X_0$ corresponds to a point of $V_{g-1, f-1} \cap
\overline{\CH}_g$. The choice of the nodal point $P_0$ is
equivalent to a choice of two distinct points of $X_0$ which are
exchanged by the hyperelliptic involution. So $\dime(C_0 \cap
\Delta_0) = \dime(V_{g-1, f-1} \cap \overline{\CH}_g) + 1$.

Furthermore, ${\rm codim}(C_0) \geq {\rm codim} (C_0 \cap \Delta_0) - {\rm codim}(\Delta_0)$
in $\overline{\CH}_g$.
A calculation shows that the codimension of $C_0$ in $\overline{\CH}_g$
is at least the codimension of $V_{g-1,f-1} \cap \overline{\CH}_{g-1}$ in $\overline{\CH}_{g-1}$.
Repeating this calculation, we see that the codimension of $C_0$ in $\overline{\CH}_g$
is at least the codimension of $V_{g-f,0} \cap \overline{\CH}_{g-f}$ in $\overline{\CH}_{g-f}$,
which by Proposition \ref{Pprk0} is $g-f$.
It follows that $V_{g,f} \cap \overline{\CH}_g$ is pure of codimension $g-f$ in $\overline{\CH}_g$.
\end{proof}

This is the main result in the paper on the $p$-rank of hyperelliptic curves.

\begin{theorem}\label{Thyperprk}
For all $g \geq 1$ and all $0 \leq f \leq g$, the locus $V_{g,f} \cap
\CH_g$ is non-empty of dimension $g-1+f$. In particular, there
exists a smooth hyperelliptic curve of genus $g$ and $p$-rank $f$.
\end{theorem}

\begin{proof}
By \cite[Proposition 2.7]{FVdG}, there exists a smooth
hyperelliptic curve $X$ of genus $g$ and $p$-rank equal to zero for all $g \geq 1$.
For $0 \leq f \leq g$, let $C_f$ be the component of $V_{g,f} \cap \overline{\CH}_g$ containing $X$.
By Proposition \ref{Phyperprk}, $C_f$ has codimension $g-f$ in $\overline{\CH}_g$.
It follows that $C_f \cap \CH_g$ has dimension $g-1+f$
since $C_f$ is not contained in the boundary of $\overline{\CH_g}$.
Now $C_{f-1}$ has codimension only $g-f+1$ in $\overline{\CH}_g$.
So the generic point of $C_f$ is a smooth hyperelliptic curve with $p$-rank exactly $f$.
\end{proof}

We now turn to the question of constructing Jacobians of curves with large $a$-number.
To do this, we first analyze the Jacobians of fibre products of hyperelliptic curves in Section \ref{S2}
and then analyze the geometry of the branch points of non-ordinary hyperelliptic curves in Section \ref{S3}.
Unless specified otherwise, the results in the next two sections are also valid in characteristic $0$
(but not in characteristic $2$).

\section{Fibre products of hyperelliptic curves} \label{S2}

Let $G$ be an elementary abelian $2$-group of order $2^n$.
In this section, we describe $G$-Galois covers $\phi:X \to \PP_k$
where $X$ is a smooth projective $k$-curve of genus $g$.
For such a cover $\phi$, we show that the Jacobian of $X$
decomposes into $2^n-1$ factors which are Jacobians as well.
We study some geometric properties of the Hurwitz space $H_{g,n}$ which parametrizes
isomorphism classes of such covers $\phi$.

\subsection{The moduli space $H_{g,n}$}

We first recall a result about
the coarse moduli space parametrizing isomorphism classes  of
$G$-Galois covers $\phi:X \to \PP_k$ where $X$ is a smooth projective $k$-curve of
genus $g$. This description is related to the theory of Hurwitz schemes and
gives a framework to describe these covers.  In particular, this framework allows one to
consider families of such covers with varying branch locus, to lift such a
cover from characteristic $p$ to characteristic 0, or to study the locus in
$\CM_g$ of curves with a certain type of action by $G$.

To be precise, let $F_{g,n}$ be the contravariant functor which associates
to any $k$-scheme $\Omega$ the set of isomorphism classes of $(\ZZ/2)^n$-Galois covers
$\phi_\Omega:X_\Omega \to \PP_\Omega$ where $X$ is a flat $\Omega$-curve
whose fibres are smooth projective curves of genus $g$
and where the branch locus $B$ of $\phi_\Omega$ is a simple horizontal divisor.
In other words, the branch locus consists of $\Omega$-points of $\PP_\Omega$ which do not intersect.
Since each inertia group is a cyclic group of order $2$, the Riemann-Hurwitz formula implies
$g=2^{n-2}|B|-2^n+1$.
The following facts about the Hurwitz scheme which coarsely represents this functor
are well-known over the complex numbers.

\begin{lemma} \label{Labstract}
\begin{description}
\item{i)} There exists a coarse moduli space $H_{g,n}$ for the functor $F_{g,n}$
which is of finite type over $\ZZ[1/2]$.
\item{ii)} There is a natural morphism $\tau: H_{g,n} \to \CM_g$
whose fibres have dimension three.
\item{iii)} There is a natural morphism $\beta:H_{g,n} \to {\mathbb P}^{|B|}$
which is proper and \'etale over the image.
\end{description}
\end{lemma}

\begin{proof} See \cite[Chapter 10]{V} for the construction of $H_{g,n}$  and
the morphisms $\tau$ and $\beta$ over $\CC$.   The corresponding statements
over $\ZZ[1/2]$ follow from \cite[Theorem 4]{W:thesis}.

We recall some of the details about the morphisms $\tau$ and
$\beta$.  The morphism $\tau$ associates to any $\Omega$-point of
$H_{g,n}$ the  isomorphism class of $X_\Omega$, where $\phi_\Omega:X_\Omega \to
\PP_\Omega$  is the corresponding cover of $\Omega$-curves. The fibres have
dimension three since $X_\Omega$ is isotrivial if and only if after an \'etale base change
from $\Omega$ to $\Omega'$ there is a projective linear transformation $\rho$ such that $\rho
\phi_{\Omega'}$ is constant,  \cite[Lemma 2.1.2]{Pr:deg}.

The morphism $\beta$ associates to any $\Omega$-point of $H_{g,n}$ the
$\Omega$-point of the configuration space  ${\mathbb P}^{|B|}$ determined by
the branch locus of the associated cover.  More  specifically, $\beta$
associates to any cover $\phi_\Omega:X_\Omega \to \PP_\Omega$ the $\Omega$-point
$[a_0: \ldots :a_{|B|}]$ of ${\mathbb P}^{|B|}$ where $a_i$ are the
coefficients of the polynomial  whose roots are the branch points of
$\phi_\Omega$. Note that the $k$-points of the image of $\beta$ correspond to
polynomials with no multiple roots.
\end{proof}

We denote by $\CH_{g,n}$ the image $\tau(H_{g,n})$ in $\CM_g$.
Given a smooth connected $k$-curve $X$, then $X$ corresponds to a
point of $\CH_{g,n}$  if and only if there exists a subgroup $G
\subset \Aut(X)$ with quotient $X/G \simeq \PP$. Note that
$\CH_{g,1}$ is simply the locus $\CH_g$ of hyperelliptic curves in
$\CM_g$.

It is often more useful to describe the branch locus of $\phi_\Omega$ directly as
an  $\Omega$-point of $(\PP)^{|B|}$.  This can be done by considering an
ordering of the branch points of $\phi_\Omega$.   The branch locus of a cover
corresponding to a $k$-point of $H_{g,n}$ can be any $k$-point of $(\PP)^{|B|}
-\Delta$ where $\Delta$ is the weak diagonal consisting of points having at
least two equal coordinates.   In particular, for any $\Omega$-point $(b_1,
\ldots, b_{2g+2})$  of $(\PP)^{2g+2} -\Delta$ there is a unique hyperelliptic
cover $\phi_{\Omega}: X_\Omega \to \PP_\Omega$ branched at $\{b_1, \ldots,
b_{2g+2}\}$.  Also the curve $X_\Omega$ has genus $g$.

\subsection{The fibres of $H_{g,n}$}

We now describe some properties of a $G$-Galois cover $\phi: X \to \PP$ corresponding
to a point of $H_{g,n}$.
In fact, the cover $\phi$ arises as the fibre product of $n$ hyperelliptic
covers which satisfy a strong disjointness condition on their branch loci.

Consider an isomorphism $\iota:(\ZZ/2)^n \simeq G$.  For $i\in \{1,\ldots,
n\}$,  this isomorphism determines a natural element $s_i$ of order $2$ in
$G$. Let $H_i \simeq (\ZZ/2)^{n-1}$ be the subgroup generated by all $s_j$ for
$j \not =i$.     Suppose for $i \in \{1, \ldots, n\}$ that $B_i$ is a
non-empty finite subset of $\PP$ of even cardinality. For any non-empty $S
\subset \{1, \ldots, n\}$, denote by  $B_S$ the set of all $b \in \PP$ such
that $b \in B_i$ for an odd number of $i \in S$ and denote by $C_S \to \PP$
the hyperelliptic cover branched at $B_S$.  Finally, let $H_S$ be the subgroup
of $G$ consisting of all elements $\sum_{i=1}^n a_i s_i$ such that $\sum_{i
\in S} a_i$ is even.  Note that each $H_S$ is non-canonically isomorphic to
$(\ZZ/2)^{n-1}$. Furthermore, when $S=\{i\}$ we have $B_S=B_i$ and $H_S=H_i$ .

\begin{lemma} \label{Lquotient}
Suppose $\phi:X \to \PP$ is the normalized fibre product over $\PP$ of
$n$ smooth hyperelliptic covers $C_i \to \PP$ with branch loci $B_i$.
Then $\phi$ is a $G$-Galois cover and the quotient of $X$ by $H_S$
is the hyperelliptic cover $C_S \to \PP$ branched at $B_S$.
\end{lemma}

\begin{proof}
The cover $\phi:X \to \PP$ is a $G$-Galois cover of (possibly
disconnected) smooth curves by the definition of the fibre product.
Also by definition, $C_i \to \PP$ is the quotient of $X$ by the
subgroup $H_i$.

The branch locus $B$ of $\phi$ equals $\cup_{i=1}^n B_i$.  For $b \in B$, the inertia group
$I_b$ of $X \to \PP$ above $b$ must be cyclic; thus $I_b \simeq \ZZ/2$.
In fact, the generator $(\alpha_1, \ldots, \alpha_n)$ of $I_b$ satisfies
$\alpha_i=1$ if and only if $b \in B_i$.  To see this, note that if
$b \in B_i$, then $C_i \to \PP$ is branched at $b$ and so $I_b \not \subset H_i$;
it follows that $\alpha_i =1$ if $b \in B_i$.  On the other hand,
if $b \not \in B_i$, then $C_i \to \PP$ is unramified at $b$ and so $I_b \subset H_i$;
it follows that $\alpha_i=0$ if $b \not \in B_i$.

Since $H_S \simeq (\ZZ/2)^{n-1}$, the quotient $X/H_S \to \PP$ is hyperelliptic;
it remains to show that the branch locus of this cover is $B_S$.
For $b \in B$, the cover $X/H_S \to \PP$ is branched at $b$ if and only if
$I_b \not \subset H_S$, which is equivalent to $(\alpha_1, \ldots, \alpha_n)
\not \in H_S$.  So $X/H_S \to \PP$ is branched at $b$ if and only if
$\Sigma_{i \in S} \alpha_i$ is odd.  Now,
$\Sigma_{i \in S} \alpha_i \equiv \#\{i \in S | b \in B_i\} \bmod 2$, so
this number is odd if and only if $b \in B_i$ for an odd number of $i \in S$.
Thus $X/H_S$ is branched at $B_S$ by definition.
\end{proof}

In Section \ref{Sanumb}, we construct covers $\phi:X \to \PP$
corresponding to points of  $H_{g,n}$ for which $X$ is also
hyperelliptic.  For example, when $n=2$, suppose $\phi$ is the
normalized fibre product of two hyperelliptic covers $\phi_1$ and
$\phi_2$.  The curve $X$ will also be hyperelliptic if its quotient
$C_{1,2} = X/H_{1,2}$ is isomorphic to $\PP$. This occurs when
$g_1=g_2$ and $B_1$ and $B_2$ overlap in all but one point; or
when $g_2=g_1 +1$ and $B_1 \subset B_2$.
The other extreme is considered in \cite{St} where
Stepanov uses the fibre product of two hyperelliptic curves
whose branch loci intersect in a single point to construct Goppa codes.

We say that the collection $\{B_i \}_{i = 1}^n$ is {\it strongly
disjoint} if the following two conditions are satisfied: first,
the sets $B_S$ are distinct for all non-empty $S \subset \{1,
\ldots, n\}$; second, $B=\cup_{i=1}^n B_i$ is a simple horizontal
divisor. In other words, if $b_1, b_2 \in B$ are two
$\Omega$-points of $\PP_\Omega$  for some scheme $\Omega$, then
the second condition insures that either $b_1=b_2$ or that $b_1$
and $b_2$ do not intersect in $\PP_\Omega$.

\begin{lemma} \label{Lfibre} A cover $\phi:X \to \PP$ corresponds to a point of
$H_{g,n}$ if and only if  $X$ has genus $g$ and $\phi:X \to \PP$ is isomorphic to
the normalized fibre product over $\PP$ of $n$ smooth hyperelliptic covers $C_i
\to \PP$ whose branch loci $B_i$ form a strongly disjoint set. \end{lemma}

\begin{proof}  If $\phi:X \to \PP$ is the normalized fibre product of $n$
hyperelliptic covers with branch loci $B_i$, then it is clear that $\phi$ is a
$G$-Galois cover and $X$ is projective.   Furthermore, $C_j \to \PP$ is
disjoint from the normalized fibre product of all $C_i \to \PP$ for $i \leq j$;
otherwise, by Lemma \ref{Lquotient}, $C_j \to \PP$ would be isomorphic  to $C_S
\to \PP$ for some $\{j\} \not = S \subset \{1, \ldots, n\}$.  This would imply
$B_S=B_j$ for some $S \not = \{j\}$ which  would contradict the fact that
$\{B_i\}$
form a strongly disjoint set. Since these covers are disjoint over $\PP$, it
follows that $X$ is connected.   Also $X$ is a smooth relative curve since
$B=\cup_{i=1}^n B_i$ is a simple horizontal divisor. By hypothesis, $X$ has
genus $g$ and so $\phi$ corresponds to a point of $H_{g,n}$.

Conversely, if $\phi:X \to \PP$ corresponds to a point of $H_{g,n}$, then $X$ has genus $g$ by
definition.  Consider the quotients $C_i \to \PP_k$ of $\phi$
by the subgroups $H_i$ of $G$ for $i=1, \ldots, n$.
These covers are clearly smooth and hyperelliptic.   By the universal property of fibre
products, there is a morphism from $X$ to the normalized fibre product of the
covers $C_i \to \PP$.   This morphism must be an isomorphism since both $X$ and
the normalized fibre product have degree $2^{n}$ over $\PP$. Also, $\phi: X \to \PP$
dominates the fibre product of any two of the quotients $C_S \to \PP$ with
branch locus $B_S$, by Lemma \ref{Lquotient}.    Since $X$ is connected, these
quotients $C_S \to \PP$ must all be disjoint; in other words, the sets $B_S$
must all be distinct.
Also, $\cup_{i=1}^n B_i$ is the branch locus $B$ of $\phi$; by definition, $B$ is
a simple horizontal divisor.  Thus $\{B_i\}$ form a strongly disjoint set.
\end{proof}

\begin{corollary} \label{Cdim} For $n \geq 2$,
the locus $\CH_{g,n}$ has dimension $(g+2^n-1)/2^{n-2}-3$ if $g
\equiv 1 \bmod 2^{n-2}$ and is empty otherwise. In particular, the
dimension of the locus $\CH_{g,2}$ is $g$.
\end{corollary}

\begin{proof} The dimension of $H_{g,n}$ is equal to the dimension of
$(\PP)^{|B|}$, namely the number of branch points $|B|$ of the corresponding
covers.   By the Riemann-Hurwitz formula, $|B| = (g+2^n-1)/2^{n-2}$.  By Lemma
\ref{Labstract}, the dimension of $\CH_{g,n}$ is three less than the dimension
of $H_{g,n}$, which simplifies to $g$ when $n=2$. \end{proof}

\subsection{Decomposition of the Jacobian} We will now describe the isogeny
class of the Jacobian for any curve $X$  for which there exists a cover $\phi:X
\to \PP_k$ corresponding to a $k$-point of $H_{g,n}$.  For $i \in \{1, \ldots, n\}$,
suppose $\phi_i:C_i \to \PP_k$ is a smooth hyperelliptic cover with branch locus
$B_i$.  Suppose $\{B_i\}_{i=1}^n$ form a strongly disjoint set and let $B=\cup_{i=1}^n
B_i$.

\begin{proposition} \label{Pfibre} Suppose $\phi:X \to \PP_k$ is the
normalization of the  fibre product of $\phi_i$ for $i =1, \ldots, n$. Then
$\Jac(X)$ is isogenous to $\prod(\Jac(C_S))$ where the product is taken over
all non-empty $S \subset \{1, \ldots, n\}$. \end{proposition}

\begin{proof}  Note that $X/H_S$ is the hyperelliptic curve $C_S$ by Lemma
\ref{Lquotient}.  Thus the result follows directly from \cite[Theorem C]{KR} if
$\gen(X)= \Sigma_S \gen(C_S)$.  By the Riemann-Hurwitz formula,
$\gen(C_S)=-1+|B_S|/2$.  Since $B=\cup_{i=1}^n B_i$ is the branch locus of $X
\to \PP_k$, it follows  that $\gen(X)=2^{n-2}|B|-2^n+1$.  The proof follows by
showing that $\Sigma_S |B_S| =2^{n-1}|B|$ by the inclusion-exclusion
principle. \end{proof}

The isogeny between $\Jac(X)$ and $\prod(\Jac(C_S))$ is not sufficient to study the $a$-number of $X$
since the $a$-number is not an isogeny invariant.  For this reason, we now generalize Proposition \ref{Pfibre}
by showing that the de Rham cohomology group $H^1_{{\rm dR}}(X)$ also decomposes.
Equivalently, one could work with the crystalline cohomology group $H^1_{\crys}(X)$ evaluated at $k$,
\cite[1.3.6]{Icrys}.
We thank Kani \cite{K} for helping us with the proof of
Proposition \ref{Pcohom}.   Let $N=2^n-1$.

\begin{proposition} \label{Pcohom} Suppose ${\rm char}(k) \not = 2$.
Then $H^1_{{{\rm dR}}}(X)$ is isomorphic to $\oplus_S H^1_{{{\rm
dR}}}(C_S)$ as $k[G]$-modules, where the sum is taken over all
non-empty $S \subset \{1, \ldots, n\}$.
\end{proposition}

\begin{proof}
Since ${\rm char}(k) \not = 2$, there exists an idempotent $\epsilon_S$ corresponding to the subgroup $H_S$
in the group ring $k[G]$ for every nonempty subset $S \subset \{1,\ldots,n\}$.
Namely, $\epsilon_S = \sum h/2^{n-1}$, where the sum ranges over all $h \in H_S$.
Let $\epsilon_G$ be the idempotent $\sum h/2^n$, where the sum ranges
over all $h \in G$. By Lemma \ref{Lquotient}, $C_S$ is the
quotient of $X$ by $H_S$, so  $H^1_{{\rm dR}}(C_S) \cong
(H^1_{{\rm dR}}(X))^{H_S} \cong \epsilon_S H^1_{{\rm dR}}(X)$.
Furthermore, note that $0 = H^1_{{\rm dR}}(\PP_k) \cong (H^1_{{\rm
dR}}(X))^G \cong \epsilon_G H^1_{{\rm dR}}(X)$ and therefore
that $\epsilon_Gx=0$ for all $x$.

If $S$ and $T$ are distinct subsets then $\epsilon_S \epsilon_T =
2^{2-2n} \sum h_sh_t$ where the sum ranges over all $h_s \in H_S$
and $h_t \in H_T$. For each $g \in G$, we see that $gh_s^{-1} \in
H_T$ for half of the values of $h_s \in H_S$. So $g$ appears
$2^{n-2}$ times in $\sum h_sh_t$. Thus,
$2^{2n-2}\epsilon_S\epsilon_T= 2^{n-2}\sum_{g \in G}g$ and we
obtain that $\epsilon_S\epsilon_T=\epsilon_G$.  Similarly, one can
show for all subsets S that $\epsilon_S \epsilon_S = \epsilon_S$
and $\epsilon_S \epsilon_G = \epsilon_G$.

We construct an explicit homomorphism $\gamma$ from $\oplus_S
H^1_{{\rm dR}} (C_S)$ to  $H^1_{{\rm dR}}(X)$:
$$\gamma(x_1,x_2,\ldots,x_N) = \sum_{i=1}^N x_i.$$
If $\psi$ is the homomorphism from $H^1_{{\rm dR}}(X)$ to
$\oplus_S H^1_{{\rm dR}}(C_S)$ given by $$\psi(y) = (\epsilon_1
y,\epsilon_2 y,\ldots, \epsilon_N y)$$ then one can check that
$\psi \circ \gamma = \gamma \circ \psi ={\rm Id}$. Thus $\gamma$
is an isomorphism of $k$-vector spaces.  In
fact, $\gamma$ is a $k[G]$-module isomorphism since every $g \in G$
commutes with $\epsilon_S$ and thus with $\gamma$.
\end{proof}

The following corollary will be used throughout the remainder of the paper.

\begin{corollary}\label{Leotype} Suppose ${\rm char}(k) > 2$. There is an isomorphism
between $\Jac(X)[p]$ and $\prod_S(\Jac(C_S)[p])$ as group schemes where the product is taken over all
non-empty $S \subset \{1, \ldots, n\}$.  In particular, $\Jac(X)$ and $\prod_S(\Jac(C_S))$ have the
same $p$-rank and $a$-number.
\end{corollary}

\begin{proof}  By Proposition \ref{Pcohom}, there is an isomorphism of $k$-vector spaces
between $H^1_{{\rm dR}}(X)$ and $\oplus_S H^1_{{\rm dR}}(C_S)$.
By the functoriality of the Frobenius and Verschiebung morphisms, $F$ and $V$ commute with the action of
$g \in G$ and thus with the idempotents $\epsilon_S$.  It follows that
$H^1_{{\rm dR}}(X)$ and $\oplus_S H^1_{{\rm dR}}(C_S)$ are naturally isomorphic as $k[V,F]$-modules.
Since $X$ and $C_S$ are smooth curves, \cite[3.11.2]{Il} implies that
$H^1_{{\rm dR}}(\Jac(X))$ and $\oplus_S H^1_{{\rm dR}}(\Jac(C_S))$ are isomorphic as $k[V,F]$-modules.
By \cite[5.11]{Oda}, $H^1_{{\rm dR}}(\Jac(X))$
is canonically isomorphic to the contravariant Dieudonn\'e module associated to $\Jac(X)[p]$.
Likewise, $H^1_{{\rm dR}}(\Jac(C_S))$)
is canonically isomorphic to the contravariant Dieudonn\'e module associated to $\Jac(C_S)[p]$.
So the Dieudonn\'e module of $\Jac(X)[p]$ is isomorphic to the direct sum of the Dieudonn\'e modules
of $\Jac(C_S)[p]$.
It follows, from the equivalence of categories between finite commutative group schemes over $k$
and their contravariant Dieudonn\'e modules,
that the group schemes $\Jac(X)[p]$ and $\prod_S(\Jac(C_S)[p])$ are isomorphic.
\end{proof}

\section{Configurations of non-ordinary hyperelliptic curves} \label{S3}

The results in this section will be used to find curves $X$
having interesting $p$-power torsion, as measured in terms of invariants such
as the $p$-rank and $a$-number. Corollary \ref{Leotype}
shows that when a cover $\phi:X \to \PP$
corresponds to a point of $H_{g,n}$ then such invariants for $X$ can be
determined by the corresponding invariants of its $\ZZ/2$-quotients.
Since these quotients are all hyperelliptic, one can apply results of Yui
\cite{Y}. The main difficulty is to control the $p$-torsion of all of the
curves $C_S$ in terms of the $p$-torsion of the curves $C_i$.

Let $C$ be a smooth hyperelliptic curve of genus $g$ defined over
an algebraically closed field $k$ of characteristic $p > 2$. Recall
that $C$ admits a $\ZZ/2$-Galois cover $\phi_1:C \to \PP_k$ with $2g+2$
distinct branch points.  Without loss of generality, we suppose
$\phi_1$ is branched at $\infty$ and choose an equation for this cover
of the form $y^2=f(x)$, where $f(x)$ is a polynomial of degree
$2g+1$. We denote the roots of $f(x)$ by
$\{\lambda_1,\ldots,\lambda_{2g+1}\}$.

Denote by $c_r$ the coefficient of $x^r$ in the expansion of
$f(x)^{(p-1)/2}$.  Then
\begin{equation}\label{expck} c_r = (-1)^{r-(p-1)/2}\sum \binom{(p-1)/2}{a_1}
\ldots\binom{(p-1)/2}{a_{2g+1}}\lambda_1^{a_1}\ldots\lambda_{2g+1}^{a_{2g+1}}
\end{equation} \noindent where the sum ranges over all $2g+1$-tuples
$(a_1,\ldots, a_{2g+1})$ of integers such that $0 \leq a_i \leq (p-1)/2$ for
all $i$ and $\sum a_i = (2g+1)(p-1)/2-r$.  Note that $c_r$ can be viewed as a
polynomial in $k[\lambda_1, \ldots, \lambda_{2g+1}]$ which is homogeneous of
degree $(2g+1)(p-1)/2 - r$ and which is of degree $(p-1)/2$ in each variable.

Let $A_g$ be the $g \times g$ matrix whose $ij$th entry is
$c_{ip-j}$. The determinant of $A_g$ defines a polynomial in
$k[\lambda_1, \ldots, \lambda_{2g+1}]$ which we denote by
$\Det_g(\lambda_1,\ldots,\lambda_{2g+1}) =
\Det_g(\vec{\lambda}_{2g+1})$. This polynomial is of degree at
most $g(p-1)/2$ in each $\lambda_i$ and is homogeneous of total
degree $g^2(p-1)/2$.
It is invariant under the action of $S_{2g+1}$ on the variables $\lambda_i$.
We denote by $D_g \subset (\Aa_k)^{2g+1}$
the hypersurface of points $\vec{\lambda}_{2g+1} =(\lambda_1,
\ldots, \lambda_{2g+1})$ for  which
$\Det_g(\vec{\lambda}_{2g+1})=0$.

In \cite{Y}, Yui gives the following characterization of non-ordinary hyperelliptic curves.
Recall that $\Delta$ is the weak diagonal consisting of points with at
least two equal coordinates.

\begin{theorem} ({\rm Yui \cite{Y}}) Suppose $C$ is a smooth hyperelliptic curve of genus $g$.
Then $C$ is non-ordinary if and only if there is a $\ZZ/2$-Galois cover
$\phi: C \to \PP_k$ branched at $\infty$ and at $2g+1$ distinct points
$\lambda_i \in \Aa_k$  such that
$\vec{\lambda}_{2g+1} \in D_g$.
\end{theorem}

We now find some results on the geometry of the hypersurface $D_g$
which will be used in Sections \ref{Sanumb} and \ref{Sptor} to
construct curves in $\CH_{g,n}$ whose $p$-torsion has prescribed
invariants.  In Lemma \ref{Ldegn} and Lemma \ref{Lrtsoffdiag},
we show that $\Det_g(\vec{\lambda}_{2g+1})$ is
generically a polynomial of degree $d=g(p-1)/2$ in the variable $\lambda_{2g+1}$
whose roots are not contained in $\{\lambda_1, \ldots, \lambda_{2g}\}$.
We expect for a generic choice of $\lambda_1,\ldots,\lambda_{2g}$ that
this polynomial will have $d$ distinct roots.
Showing this seems to be related to the question of whether the hyperelliptic locus
is transversal (in the strict geometric sense) to the locus $V_{g,g-1}$ of nonordinary curves.
So in Proposition \ref{Pexactlyd}, we instead prove the weaker statement
that this polynomial has at least $(p-1)/2$ distinct roots.

\begin{lemma} \label{Ldegn}
The determinant $\Det_g(\vec{\lambda}_{2g+1})$ is a polynomial of
degree $d=g(p-1)/2$ in the variable $\lambda_{2g+1}$.
\end{lemma}

\begin{proof}  As we observed above, the degree of $\Det_g(\vec{\lambda}_{2g+1})$ in
$\lambda_{2g+1}$ is at most $d$. We claim that the coefficient of
$\lambda_{2g+1}^d$ is a non-zero polynomial in $k[\lambda_1, \ldots,
\lambda_{2g}]$.  In particular, one term of this polynomial is
$(-1)^{g(p-1)/2} \lambda_{2g+1}^d\prod_{i=1}^{2g}
\lambda_i^{(g-\lceil i/2 \rceil)(p-1)/2}$.

To see this, we note first from Equation \ref{expck} that the
total degree of $c_{gp-j}$ is $(2g+1)(p-1)/2 - (gp-j) = (p-1)/2 +
(j-g)$. So if $j < g$ then $\lambda_{2g+1}^{(p-1)/2}$ cannot
appear in $c_{gp-j}$. Furthermore, the coefficient of
$\lambda_{2g+1}^{(p-1)/2}$ in $c_{gp-g}$ is exactly
$(-1)^{(p-1)/2}$. Because the degree of $\lambda_{2g+1}$ in $c_r$
is at most $(p-1)/2$ for all $r$, a monomial in
$\Det_g(\vec{\lambda}_{2g+1})$  will be divisible by
$\lambda_{2g+1}^d$ only if it is the product of matrix entries
which are each divisible by $\lambda_{2g+1}^{(p-1)/2}$.  Thus
$c_{gp-g}$ is the only entry in the bottom row of $A_g$ which
contributes to the  terms of $\Det_g(\vec{\lambda}_{2g+1})$ which
are divisible by $\lambda_{2g+1}^d$.

Similarly, in the penultimate row of $A_g$, the total degree of
$c_{(g-1)p-j}$ will be $3(p-1)/2 + (j-(g-1))$.  Therefore, if $j <
g-1$ then $(\lambda_1\lambda_2\lambda_{2g+1})^{(p-1)/2}$ cannot
divide $c_{(g-1)p-j}$. Because the degree of $\lambda_1$ for all
$c_r$ in $A_g$ is at most $(p-1)/2$, only the last two
entries of the penultimate row contribute to the terms of
$\Det_g(\vec{\lambda}_{2g+1})$ which are divisible by
$\lambda_1^{(g-1)(p-1)/2}$.  Also the coefficient of
$(\lambda_1\lambda_2\lambda_{2g+1})^{(p-1)/2}$ in
$c_{(g-1)p-(g-1)}$ is $(-1)^{(p-1)/2}$.

Continuing, we see that only terms which are on or above the
diagonal can contribute to the desired term of
$\Det_g(\vec{\lambda}_{2g+1})$. It follows that the only term of
$\Det_g(\vec{\lambda}_{2g+1})$ which involves the monomial
$\lambda_{2g+1}^d\prod_{i=1}^{2g} \lambda_i^{(g-\lceil i/2
\rceil)(p-1)/2}$ comes from the product of elements of the
diagonal. The coefficient of this monomial is the product of $g$
coefficients which each equal $(-1)^{(p-1)/2}$, so it is equal to
$(-1)^{g(p-1)/2}$.
\end{proof}

\begin{lemma} \label{Lrtsoffdiag}
The image of $\Det_g(\vec{\lambda}_{2g+1})$ in $k[\lambda_1,
\ldots, \lambda_{2g+1}]/(\lambda_{2g+1}-\lambda_1)$ is
non-constant.
\end{lemma}

\begin{proof}
The proof is similar to that of Lemma \ref{Ldegn}.
It is sufficient to show that at least one of the coefficients of
$\Det_g(\lambda_1,\ldots, \lambda_{2g},\lambda_1)$ is non-zero.
The coefficient of the monomial $\lambda_{1}^{g(p-1)/2}\prod_{i=1}^{2g}
\lambda_i^{(g-\lceil i/2 \rceil)(p-1)/2}$ is $2(-1)^{g(p-1)/2}$ as
this monomial appears exactly twice as the product of terms in the diagonal of the
Hasse-Witt matrix and does not appear again in the expansion of the determinant.
\end{proof}

Suppose exactly two branch points of a smooth hyperelliptic cover specialize together.
The resulting curve is singular and consists of a hyperelliptic curve $C'$
of genus $g-1$ self-intersecting in a point.
The geometric interpretation of the next lemma is that
this singular curve will be ordinary if and only if $C'$ is ordinary.

\begin{lemma}\label{Lzeros}
$\Det_g(\lambda_1,\ldots,\lambda_{2g-1},0,0) =
(- \lambda_1 \cdots \lambda_{2g-1})^{(p-1)/2}\Det_{g-1}(\lambda_1,\ldots,
\lambda_{2g-1})$.
\end{lemma}

\begin{proof}
Suppose $\lambda_{2g}=\lambda_{2g+1}=0$.
Then the only nonzero terms in the sum defining $c_r$ are
those where $a_{2g}=a_{2g+1}=0$.
If $r=p-1$, then the only term in this sum that does
not vanish is the one where $a_i=(p-1)/2$ for $1 \leq i \leq 2g-1$.
Thus $c_{p-1} = (-\lambda_1 \cdots\lambda_{2g-1})^{(p-1)/2}$.
If $r<p-1$, then all of the terms in the sum are zero, and thus $c_r=0$.
Suppose $r > p-1$ and $r=ip-j$.
Then the term $c_r$ occuring in the $i$th row and $j$th column of $A_g$ equals the term
$c_{r-(p-1)}$ occuring in the $(i-1)$st row and $(j-1)$st column of $A_{g-1}$.
By expanding the determinant along the first row, we see that
$\Det_g(\lambda_1,\ldots,\lambda_{2g-1},0,0) =c_{p-1} \Det(A_{g-1})$.
\end{proof}

For fixed $\vec{\lambda}_{2g}=(\lambda_1, \ldots, \lambda_{2g})
\subset (\Aa_k)^{2g}$, denote by $L(\vec{\lambda}_{2g})$ the line
consisting of points $(\lambda_1, \ldots, \lambda_{2g},
\lambda_{2g+1}) \subset (\Aa_k)^{2g+1}$ (where only the last
coordinate varies).
Generically, the intersection of $L(\vec{\lambda}_{2g})$ and $D_g$ consists of
$d=g(p-1)/2$ points when counted with multiplicity.
To see this, consider $\Det(\vec{\lambda}_{2g+1})$ as a polynomial in $R[\lambda_{2g+1}]$
where $R=k[\lambda_1, \ldots, \lambda_{2g}]$.
The coefficient of $\lambda_{2g+1}^d$ in $\Det(\vec{\lambda}_{2g+1})$
is non-zero in $R$ by Lemma \ref{Ldegn}.
Since $k$ is an algebraically closed field,
for any $\vec{\lambda}_{2g}=(\lambda_1, \ldots, \lambda_{2g})$
not in the Zariski closed set of $(\Aa_k)^{2g}$ defined by this coefficient,
$\Det(\vec{\lambda}_{2g+1})$ has degree $d$ and thus $d$ roots
in $k$ when counted with multiplicity.
The next proposition gives a lower bound on the number of distinct roots.

\begin{proposition}  \label{Pexactlyd} Let $U_g \subset (\Aa_k)^{2g}$ be the
set of points $(\lambda_1, \ldots, \lambda_{2g})$ for which
$L(\vec{\lambda}_{2g})$ intersects $D_g$ in at least $(p-1)/2$
non-zero distinct points of $(\Aa_k)^{2g+1} \backslash \Delta$. Then $U_g$
is a nonempty Zariski open subset of $(\Aa_k)^{2g}$.
\end{proposition}

\begin{proof}
The proof is by induction on $g$.  A result of Igusa \cite{I} states
that there are exactly $(p-1)/2$ distinct values $\lambda$ so
that the elliptic curve branched at $\{0,1,\infty, \lambda\}$ is non-ordinary.
It follows that the result is true when $g=1$.

Suppose that $U_{g-1}$ is a nonempty Zariski open subset of
$(\Aa_k)^{2g-2}$. First we show that for a generic choice of
$(\lambda_1,\ldots,\lambda_{2g})$ there are at least $(p-1)/2$
distinct choices of $\lambda_{2g+1}$ so that
$\Det_g(\lambda_1,\ldots,\lambda_{2g},\lambda_{2g+1}) = 0$. It
will suffice to construct a single choice of
$(\lambda_1,\ldots,\lambda_{2g})$ for which this result holds, as
the generic case will have at least as many distinct roots as any
specialized case.  It follows from Lemma \ref{Lzeros} that for
non-zero $\lambda_3, \ldots, \lambda_{2g}$,
 the non-zero values of $\lambda_{2g+1}$ so that
$\Det_g(0,0,\lambda_3,\ldots,\lambda_{2g},\lambda_{2g+1})=0$ and
$\Det_{g-1}(\lambda_3, \ldots, \lambda_{2g+1})=0$ are the same. By
the inductive hypothesis, for the generic choice of
$(\lambda_3,\ldots,\lambda_{2g})$ there are at least $(p-1)/2$
non-zero distinct values of $\lambda_{2g+1}$ with this property.

Next we show that generically these $(p-1)/2$ intersection points of $L(\vec{\lambda}_{2g})$ and $D_g$
are not contained in $\Delta$.
By Lemma \ref{Lrtsoffdiag} and by symmetry, for each $1 \leq i \leq 2g$, the value
$\lambda_i$ is a root of the polynomial $\Det_g(\lambda_1, \ldots, \lambda_{2g+1}) \in R[\lambda_{2g+1}]$
only when $(\lambda_1, \ldots, \lambda_{2g})$ is in a Zariski closed subset of $(\Aa)^{2g}$.
So for the generic choice of $(\lambda_1,\ldots,\lambda_{2g})$,
the root $\lambda_{2g+1}$ will not be contained in $\{\lambda_1, \ldots, \lambda_{2g}\}$.
It follows that for the generic choice of $(\lambda_1, \ldots, \lambda_{2g})$ the line
$L(\vec{\lambda}_{2g})$ intersects $D_g$ in at least $(p-1)/2$
non-zero distinct points of $(\Aa_k)^{2g+1} \backslash \Delta$.  So $U_g$
is a nonempty Zariski open subset of $(\Aa_k)^{2g}$.
\end{proof}

\begin{proposition} \label{Paddtwo} Let $U_{g} \subseteq (\Aa_k)^{2g}$
be defined as in Proposition \ref{Pexactlyd}. Then we have that
$U_{g} \cap (D_{g-1} \times \Aa_k)$ has codimension 1 in
$(\Aa_k)^{2g}$.
\end{proposition}

\begin{proof} Since $D_{g-1} \times \Aa_k$ has codimension 1 in $(\Aa)^{2g}$
and $U_g$ is open by Proposition \ref{Pexactlyd}, it is sufficient
to show that no component $V$ of $D_{g-1} \times \Aa_k$ is
contained in the complement $W_g$ of $U_g$. Note that the
complement of $U_g$ is a Zariski closed subset  defined by
equations which are each symmetric in the variables  $\lambda_1,
\ldots, \lambda_{2g}$.  On the other hand, any component $V$ of
$D_{g-1} \times \Aa_k$ is defined by equations that do not involve
$\lambda_{2g}$.   Since the ideal of $W_g$ is not contained in the
ideal of $V$, it follows that $V$ is not contained in $W_g$.
\end{proof}

\section{Curves with prescribed $a$-number} \label{Sanumb}

We now consider the $a$-number of Jacobians of curves with commuting involutions.
Recall that the $a$-number, $\dime_k \Hom(\alpha_p, \Jac(X))$, of a $k$-curve $X$
is an integer between $0$ and $g$.
Here $\alpha_p$ is the kernel of Frobenius on $\GG_a$.
A generic curve is ordinary and thus has $a$-number equal to zero. A supersingular
elliptic curve $E$ has $a$-number equal to one and in this case
there is a non-split exact sequence $0 \to \alpha_p \to E[p] \to
\alpha_p \to 0$.
There is a unique isomorphism type of group scheme for
the $p$-torsion of a supersingular elliptic curve, which we denote $M$.
In this section we construct curves $X$ so that $\Jac(X)[p]$ contains
multiple copies of the group scheme $M$ and thus has large $a$-number.

Let $T_{g,a}$ denote the sublocus of $\overline{\CM}_g$ consisting of curves of genus $g$ with
$a$-number at least $a$. The codimension of $T_{g,a}$ in $\CM_g$
is at least $a$ since $T_{g,a} \subset V_{g, g-a}$. It is not
known whether (for all $g$ and all $0 \leq a \leq g$) there exists
a curve of genus $g$ with  $a$-number equal to $a$.  The results
in this section give some evidence for a positive answer to this
question.

We note that these results can be viewed as a generalization of
\cite[Section 5]{O}.  In that paper, Oort considers curves $X$ of
genus $g=3$ with a group action by $G = (\ZZ/2)^2$ so that the
three $\ZZ/2$-quotients of $X$ are all elliptic curves. He shows
that there exist (nonhyperelliptic) curves of genus $3$ with
$a$-number $3$ for all primes $p \ge 3$ as well as hyperelliptic
supersingular curves of genus $3$ with $a$-number $3$ for all $p
\equiv 3 \bmod 4$.

\begin{lemma} \label{lord}
The generic geometric point of the hyperelliptic locus $\CH_g$ has $a$-number equal to $0$.
The non-ordinary locus has codimension one in $\CH_g$ and its generic geometric point
has $a$-number $1$ and $p$-rank $g-1$.
\end{lemma}

\begin{proof}
This is immediate from Theorem \ref{Thyperprk} and the fact that a curve with $p$-rank $g-1$ has
$a$-number $1$.
\end{proof}

The next theorem will lead immediately to Corollary \ref{Pa=n}
which is the main result in this paper on the $a$-number of curves.

\begin{theorem} \label{Ta=n}
Suppose $n \geq 2$ and $p \geq 2n+1$.
Suppose $g$ is such that $g \equiv 1 \bmod 2^{n-2}$ and $g \ge (n-1)2^{n-2}+1$.
There exists a family of smooth curves $X$ of genus $g$ of dimension at least $(g+2^n-1)/2^{n-2}-3-n$
so that $\Jac(X)[p]$ contains the group scheme $M^n$.
\end{theorem}

For the proof of Theorem \ref{Ta=n}, we will construct a fibre
product $\phi:X \to \PP$ of $n$ hyperelliptic covers $\phi_i$ so that
the disjoint union of any two of the branch loci $B_i$ will
consist of exactly two points.  It follows that the curves
$C_{i,j}$ will have genus zero.

\begin{proof}   Write $g = 1+\ell2^{n-2}$.
If $\ell \not\equiv n \bmod 2$, let $g_1=(\ell+3-n)/2$.  Note that $g_1 \ge 1$.
By Proposition \ref{Pexactlyd} and Lemma \ref{lord}, as long as $n \leq (p-1)/2$,
there is a Zariski open subset $U_{g_1}$ of $(\Aa_k)^{2g_1}$ with the following property:
there are at least $n$ choices $\eta_1, \ldots, \eta_n$ for $\lambda_{2g_1+1}$ such that
the corresponding hyperelliptic curve $C_i$ is non-ordinary.
By Theorem \ref{Thyperprk}, after replacing $U_{g_1}$ with a smaller
Zariski open subset of $(\Aa_k)^{2g_1}$,
we can further suppose that the curves $C_1, \ldots, C_n$ will have $p$-rank $g_1-1$.
Thus $\Jac(C_i)[p]$ contains $M$.

Let $\phi_i:C_i \rightarrow \PP$ for $1 \leq i \leq n$ be
the hyperelliptic $U_{g_1}$-curves corresponding to these choices.
Let $\phi:X \to \PP$ be the normalization of the fibre product of the
covers $\phi_i$. Note that $X$ is branched at $B=\{\infty, \lambda_1,
\ldots, \lambda_{2g_1}, \eta_1, \ldots, \eta_n\}$. By Proposition
\ref{Pfibre}, the genus of $X$ will be $2^{n-2}(2g_1+1+n)-2^n+1 =g$.
By Corollary \ref{Leotype}, $\Jac(X)[p]$ contains $\oplus_{i=1}^n \Jac(C_i)[p]$
which contains $\oplus_{i=1}^n M$.
The dimension of this family of curves is
$2g_1-2 = |B|-3-n$ which equals $(g+2^n-1)/2^{n-2}-3-n$.
Note that the $p$-rank of $X$ is at least $n(g_1-1)$.

Alternatively, if $\ell \equiv n \bmod 2$, let $g_1 =
(\ell+2-n)/2$ and note $g_1 \ge 1$.   By Proposition
\ref{Paddtwo}, the locus $U_{g_1+1} \cap (D_{g_1} \times \Aa_\ell)$ has
codimension  1 in $(\Aa_k)^{2g_1+2}$.
In other words, as long as $n \leq 1+(p-1)/2$, for any
$(\lambda_1 , \ldots,  \lambda_{2g_1+2})$ in a codimension 1 subset $Z$ of $(\Aa_k)^{2g_1+2},$
it is true that $(\lambda_1, \ldots, \lambda_{2g_1+1}) \in D_{g_1}$ and there are at least $n-1$
choices $\eta_i$ of $\lambda_{2g_1+3}$ with $(\lambda_1, \ldots, \lambda_{2g_1+3}) \in D_{g_1+1}$.
Let $\phi_n:C_n \to \PP$ be the hyperelliptic cover branched at $(\lambda_1, \ldots, \lambda_{2g_1+1})$.
For $1 \leq i \leq n-1$, let $\phi_i:C_i \to \PP$ be the hyperelliptic cover
branched at $(\lambda_1, \ldots, \lambda_{2g_1+2}, \eta_i)$.
Then $C_n$ has genus $g_1$ and $C_i$ has genus $g_1+1$ for $1 \leq i \leq n-1$.
By Theorem \ref{Thyperprk}, after restricting to a Zariski open subset of $Z$,
we can further suppose that $C_n$ (resp.\ $C_i$) has $p$-rank $g_1-1$ (resp.\ $g_1$).
Thus $\Jac(C_i)[p]$ contains $M$ for $1 \leq i \leq n$.

Let $\phi:X \to \PP$ be the normalization of the fibre product of $\phi_i$ for $1 \leq i \leq n$.
Note that $\phi$ is branched at $B=\{\infty, \lambda_1, \ldots, \lambda_{2g_1+2},
\eta_1, \ldots, \eta_{n-1}\}$. As above, $X$ has genus $2^{n-2}(2g_1+2+n)-2^n+1 = g$
and $\Jac(X)[p]$ contains $M^n$.
By Proposition \ref{Paddtwo}, the locus $Z$ has dimension $2g_1+1$.
The corresponding family of curves has
dimension $2g_1-1=|B|-3-n$ which again equals $(g+2^n-1)/2^{n-2}-3-n$.
Note that the $p$-rank of $X$ is at least $ng_1-1$.
 \end{proof}

\begin{corollary} \label{Pa=n}
Suppose $n \geq 2$ and $p \geq 2n+1$.
Suppose $g$ is such that $\CH_{g,n}$ is non-empty of dimension at least $n+1$.
Then the intersection $\CH_{g,n} \cap T_{g,n}$ has codimension at most $n$ in $\CH_{g,n}$.
In particular, there exists a smooth curve of genus $g$ with $a$-number at least $n$.
\end{corollary}

\begin{proof} By Corollary \ref{Cdim},
the condition that $\CH_{g,n}$ is non-empty is equivalent to $g \equiv 1 \bmod 2^{n-2}$ and
the condition that $\CH_{g,n}$ has dimension at least $n+1$ is equivalent to $g \ge (n-1)2^{n-2}+1$.
The family constructed in Theorem \ref{Ta=n} has dimension $(g+2^n-1)/2^{n-2}-3-n$ and thus codimension $n$ in $\CH_{g,n}$.
For any fibre $X$ in this family, $\Jac(X)[p]$ contains $M^n$ and so $X$ has $a$-number at least $n$.
So this family is contained in $\CH_{g,n} \cap T_{g,n}$.
\end{proof}

When $n=2$ or $n=3$, then the curves found in Theorem \ref{Ta=n}
are in fact hyperelliptic.

\begin{corollary} \label{Tanumb2} Suppose $g \geq 2$ and $p \geq 5$.
There exists a $(g-2)$-dimensional family of smooth
hyperelliptic curves of genus $g$ whose fibres have
$a$-number $2$ and $p$-rank $g-2$.  \end{corollary}

The family in Corollary \ref{Tanumb2} has codimension $2$ in $\CH_{g,2}$.

\begin{proof}
This follows immediately from Theorem \ref{Pa=n} once we show that
the curve $X$ is hyperelliptic when $n=2$.
If $g$ is even, note that the branch
loci of $\phi_1$ and $\phi_2$ differ in only one point. The third
quotient $C_{1,2}$ of $X$ by $\ZZ/2$ is branched at only two
points $\eta_1$ and $\eta_2$.  So the cover $X \to C_{1,2}$ is
hyperelliptic. Likewise, if $g$ is odd, then the third quotient
$C_{1,2}$ of $X$ by $\ZZ/2$ is branched at only two points
$\lambda_{2g_1+2}$ and $\eta_1$  so the cover $X \to C_{1,2}$ is
hyperelliptic.  In both cases, $\Jac(X)[p] \simeq \Jac(C_1)[p] \oplus \Jac(C_2)[p]$
and so the fibres of $X$ have $a$-number 2 and $p$-rank $g-2$.
\end{proof}

\begin{corollary} \label{Canum3} Suppose $g \geq 5$ is odd and $p \geq 7$.
There exists a $(g-5)/2$-dimensional family of smooth
hyperelliptic curves $X$ of genus $g$ so that $\Jac(X)[p]$
contains $M^3$ and thus has $a$-number at least $3$.
\end{corollary}

The family in Corollary \ref{Canum3} has codimension $3$ in $\CH_{g,3}$.

\begin{proof} It is sufficient to show that the fibres of the family
constructed in Theorem \ref{Pa=n} are hyperelliptic when $n=3$.  In both
cases of the construction, if $S=\{1,2\}$, $\{1,3\}$, or
$\{2,3\}$, then the quotient $C_S \to \PP$ of $X$ by $H_S \simeq
(\ZZ/2)^2$ is branched at only two points and so $C_S$ has genus
0.  Consider the quotient $X'$ of $X$ by the subgroup  $H' \simeq
\ZZ/2$ generated by $h=(1,1,1)$.  Note that $X'$ dominates
$C_S$ if  $S=\{1,2\}$, $\{1,3\}$, or $\{2,3\}$, since $h \in H_S$.
It follows that $X'$ has genus zero since $X' \to \PP$ is a
$(\ZZ/2)^2$-cover  having three $\ZZ/2$-quotients of genus zero.
It follows that $X \to X'$ is hyperelliptic.
\end{proof}

\begin{remark} \label{Rexact}
One would like to strengthen Corollary \ref{Canum3} by producing curves with $a$-number exactly 3.
The difficulty is to determine the $a$-number of $C_{\{1,2,3\}}$.
For example, to construct a curve of genus $g=5$ and $a$-number exactly $3$ with this method,
one would need to guarantee that there are supersingular values $\lambda_1$, $\lambda_2$
and $\lambda_3$ so that the hyperelliptic curve of genus two
branched at $\{0,1,\infty, \lambda_1, \lambda_2, \lambda_3\}$ is ordinary.
\end{remark}

\begin{remark}
In the above results, some restriction on $p$ is unavoidable.
By Proposition $3.1$ of \cite{R}, there does not exist a hyperelliptic
curve of genus 2 and $a$-number 2 when $p=3$ or of genus 3 and $a$-number 3 when $p=3$ or $5$.
Also, there does not exist a hyperelliptic curve with $a$-number 4 when $g=4$
and $p=3,5$ or when $g=5$ and $p=3$.
\end{remark}

We now produce curves of every genus with $a$-number at least 4 using this method.
(One can also produce curves of every genus with $a$-number at least 3 and count the dimension
of these families).  The curves constructed in this way
will most likely not be hyperelliptic.  This makes it difficult
to produce a curve of every genus with every possible $a$-number using induction and fibre products.

\begin{corollary}\label{Canumb3}  Suppose $g \geq 7$ and $p \geq 5$.
There exists a curve of genus $g$ with $a$-number at least $4$. \end{corollary}

\begin{proof} If $g$ is even, let $g_1 = g/2$.  Note that $g_1-2 \geq 2$.
From Corollary \ref{Tanumb2}, there exists a hyperelliptic curve of genus $g_1-2$ and $a$-number $2$.
Consider the corresponding hyperelliptic cover $\phi_1$ branched at
$\{\lambda_1, \ldots, \lambda_{2g_1-3}, \infty\}$.
Consider a hyperelliptic cover $\phi_2$ branched at $\{\eta_1, \ldots, \eta_5, \infty\}$ which has
$a$-number $2$.
After modifying $\phi_2$ by an affine linear transformation,
one can suppose that $\{\eta_i\} \cap \{\lambda_i\}$ is empty.
The cardinality of $(B_1 \cup B_2)\backslash (B_1 \cap B_2)$ is $(2g_1-2)+6-2 = 2g_1+2$.
It follows from Proposition \ref{Pfibre} that the fiber product of $\phi_1$ and $\phi_2$
yields a curve with genus $(g_1-2)+g_1+2=g$ and $a$-number at least $4$.

If $g$ is odd, let $g_1 = (g-1)/2$.  Note that $g_1-1 \geq 2$.  By
Corollary \ref{Tanumb2}, there exists a hyperelliptic curve of
genus $g_1-1$ and $a$-number $2$. Consider the corresponding
hyperelliptic cover $\phi_1$ branched at $\{\lambda_1, \ldots,
\lambda_{2g_1-2}, 0, \infty\}$. Consider a hyperelliptic cover
$\phi_2$ branched at $\{\eta_1, \ldots, \eta_4, 0, \infty\}$ which
has $a$-number $2$. After modifying $\phi_2$ by a scalar
transformation, one can suppose that $\{\eta_i\} \cap
\{\lambda_i\}$ is empty. The cardinality of $(B_1 \cup
B_2)\backslash (B_1 \cap B_2)$ is $2g_1+6-4 = 2g_1+2$. It follows
from Proposition \ref{Pfibre} that the fiber product of $\phi_1$ and
$\phi_2$ yields a curve with genus $(g_1-1)+g_1+2=g$ and $a$-number
at least $4$.
\end{proof}

\section{Curves with prescribed $p$-torsion} \label{Sptor}

The methods of the previous sections can also be used to
construct Jacobians whose $p$-torsion contains group schemes other than $\mu_p$ or $\alpha_p$.
In this section, we illustrate this for two particular isomorphism types of group scheme,
namely the $p$-torsion of a supersingular abelian surface which is not superspecial and of
a supersingular abelian variety of dimension $3$ with $a$-number $1$.

Section \ref{S2} allows one to describe the $p$-torsion of the
Jacobian of a curve $X$ which corresponds to a point of $\CH_{g,n}$.
Specifically, Proposition \ref{Pcohom} states that $\Jac(X)[p]$ is the direct sum of $\Jac(C_S)[p]$
where $C_S$ is the quotient of $X$ by $H_S$ and $S$ ranges over
the $2^n-1$ nonempty subsets of $\{1,\ldots,n\}$.
With this method, it is only possible to construct Jacobians so that
$\Jac(X)[p]$ is decomposable into (at least two) group schemes each of which occurs
as the $p$-torsion of a hyperelliptic curve.

Via the $p$-rank, we have already considered the group scheme for the
$p$-torsion of an ordinary elliptic curve, namely $\ZZ/p \oplus \mu_p$.
Using the $a$-number, we have already studied the group scheme $M$
of the $p$-torsion of a supersingular elliptic curve.

Not many other group schemes are known to occur as the $p$-torsion of a hyperelliptic curve.
There are four possibilities of group scheme which occur among curves of genus
$2$ (which are automatically hyperelliptic). The first three
$(\ZZ_p \oplus \mu_p)^2$, $(\ZZ/p \oplus \mu_p) \oplus M$, and
$M^2$ are decomposable.  We will focus on the most interesting of
the four, namely the group scheme $N$ for the $p$-torsion of a
supersingular abelian surface which is not superspecial. A curve
$X$ with $\Jac(X)[p] \simeq N$ has genus $2$ and is thus hyperelliptic.

By \cite[Example A.3.15]{G:book}, there is a filtration
$H_1 \subset H_2 \subset N$ where $H_1 \simeq \alpha_p$,
$H_2/H_1 \simeq \alpha_p \oplus \alpha_p$ and $N/H_2 \simeq \alpha_p$.
Moreover, the kernel $G_1$ of Frobenius and the kernel $G_2$ of Verschiebung are contained in $H_2$
and there is an exact sequence $0 \to H_1 \to G_1 \oplus G_2 \to H_2 \to 0$.

The group scheme $N$ is perhaps easier to describe in terms of its covariant Dieudonn\'e module.
Consider the non-commutative ring $E=W(k)[F,V]$ with the Frobenius automorphism $\sigma:W(k) \to W(k)$
and the relations $FV=VF=p$ and $F\lambda=\lambda^\sigma F$ and $\lambda V=V \lambda^\sigma$ for all
$\lambda \in W(k)$.
Recall that there is an equivalence of categories between finite commutative group schemes
${\mathbb G}$ over $k$ (with order $p^r$) and finite left $E$-modules $D({\mathbb G})$
(having length $r$ as a $W(k)$-module), see for example \cite[A.5]{G:book}.
By \cite[Example A.5.1-5.4]{G:book}, $D(\mu_p)=k/k(V, 1-F)$, $D(\alpha_p)=k/k(F,V)$,
and $D(N)=k/k(F^3,V^3,F^2-V^2)$.


The $p$-rank of a curve $X$ with $\Jac(X)[p] \simeq N$ is zero.
To see this, note that $\Hom(\mu_p, N)=0$ or that
$F$ and $V$ are both nilpotent on $D(N)$.
The $a$-number of a curve $X$ with $\Jac(X)[p] \simeq N$ is one.
(It is at least one since the $p$-rank is $0$
and at most one since the abelian surface is not superspecial.)
This also follows from the structure of the group scheme or from
the fact that $N[F] \cap N[V] =H_1 \simeq \alpha_p$.

\begin{lemma} \label{Lgenus2}
There is a one-dimensional family of smooth curves $X$ of genus two
with $\Jac(X)[p] \simeq N$.
\end{lemma}

\begin{proof}
The dimension in $\CA_2$ of supersingular (resp.\ superspecial) abelian surfaces
is one (resp.\ zero).  It follows that the locus of abelian surfaces with $p$-torsion
$N$ is exactly one.  The generic point of this one-dimensional family must be in the image of the Torelli
morphism since $\overline{\CM}_2$ and $\overline{\CA}_2$ have the same dimension.
So there is a one-dimensional family of curves of genus two with $p$-rank $0$ and $a$-number $1$.
The fibres of this family are all smooth since
the family cannot intersect either of the boundary components $\Delta_0$ or $\Delta_1$.
\end{proof}

\begin{lemma} \label{Cgenus3}
There exists a one-dimensional family of smooth hyperelliptic curves $X$ of genus $3$
with $\Jac(X)[p] \simeq N \oplus (\ZZ/p \oplus \mu_p)$.
\end{lemma}

\begin{proof}
By Lemma \ref{Lgenus2}, there is a one-dimensional family of smooth curves $X$ of genus two
with $\Jac(X)[p] \simeq N$.  This yields a family of hyperelliptic covers of $\PP$
branched at six points.  For some subset of four of these points,
the family of elliptic curves branched at these points must have varying $j$-invariant
and so its fibres are generically ordinary.  The fibre product of these two families of covers yields
a family of smooth hyperelliptic curves of genus $3$
with $p$-torsion $N \oplus (\ZZ/p \oplus \mu_p)$ by Corollary \ref{Leotype}.
\end{proof}

The following proposition will be used to generalize Lemma \ref{Cgenus3} for $g \geq 4$.

\begin{proposition} \label{TaddL}  Suppose there exists an $r$-dimensional
family of smooth hyperelliptic curves $C$ of genus $g'$ with
$\Jac(C)[p] \simeq {\mathbb G}$  for some group scheme ${\mathbb G}$.
Suppose $s \geq 1$ and $g=2g'-1+s$.  Then there exists an
$(r+s)$-dimensional family of smooth curves $X$ in $\CH_{g,2}$ so
that $\Jac(X)[p]$ contains ${\mathbb G}$.
\end{proposition}

\begin{proof} For each curve $C$ in the original family with
$\Jac(C)[p] \simeq {\mathbb G}$ and branch locus $B_0=\{\lambda_1,
\ldots, \lambda_{2g'+2}\}$, we will construct an $s$-dimensional
family of smooth curves $X$ so that $\Jac(X)[p]$ contains ${\mathbb G}$.
By Proposition \ref{Pfibre} and Lemma \ref{Leotype}, it will suffice to construct
hyperelliptic curves $C_1$ and $C_2$ whose branch loci $B_1$ and
$B_2$ are of even cardinality with $|B_1 \cap B_2|=s$ and $B_0= (B_1
\cup B_2) \backslash (B_1 \cap B_2)$.

If $s = 2m$ is even, then $B_1 = B_0 \cup
\{\eta_1,\ldots,\eta_{2m}\}$ and $B_2 =
\{\eta_1,\ldots,\eta_{2m}\}$ satisfy these restrictions and there
are $2m=s$ choices for the points $\eta_i$. Similarly, if $s =
2m+1$ is odd, then we can set $B_1 =\{\lambda_1,\ldots,
\lambda_{2g'+1},\eta_1,\ldots,\eta_{2m+1}\}$ and $B_2 =
\{\lambda_{2g'+2},\eta_1,\ldots,\eta_{2m+1}\}$ satisfy these
restrictions. There are $2m+1=s$ choices for $\eta_i$. The
Jacobian of the normalized fibre product $X$ of $C_1$ and $C_2$
contains $\Jac(C)$.
\end{proof}

This is the main result of the section.

\begin{corollary} \label{Cdim2}
Let $N$ be the $p$-torsion of a supersingular abelian surface which is not superspecial.
For all $g \geq 2$, there exists a smooth hyperelliptic curve $X$ so that
$\Jac(X)[p]$ contains $N$.
\end{corollary}

\begin{proof}
The statement will follow from induction.
Assume for all $g'$ such that $2^n \leq g' < 2^{n+1}$ that there exists a smooth hyperelliptic curve $X_{g'}$ so that
$\Jac(X_{g'})[p]$ contains $N$.
This is true for $n=1$ by Lemma \ref{Lgenus2} and Lemma \ref{Cgenus3}.
If $2^{n+1} \leq g < 2^{n+2}$, then $g=2g'$ or $g=2g'+1$ for some $g'$ such that $2^n \leq g' < 2^{n+1}$.
Using Proposition \ref{TaddL} with $s=1$ or $s=2$ allows one
to construct a curve $X_g$ of genus $g$ so that $\Jac(X_g)[p]$ contains $N$.
If $s=1$ or $s=2$ in Proposition \ref{TaddL}, then $B_2$ consists of exactly two points
so $X_{g}$ is also hyperelliptic.
\end{proof}

Similarly, one can consider the group scheme $Q$ of the $p$-torsion
of a supersingular abelian variety of dimension three with $a$-number $1$.
A curve $X$ with $\Jac(X)[p] =Q$ has $p$-rank $0$.  Also, $D(Q)=k[F,V]/k(F^4,V^4,F^3-V^3)$.
The restriction on $g$ in the next corollary could be removed if there exists a smooth hyperelliptic curve
$X$ of genus $4$ so that $\Jac(X)[p]$ contains $Q$.

\begin{corollary} \label{CQ}
Let $Q$ be the $p$-torsion of a supersingular abelian variety of dimension three with $a$-number $1$.
Suppose $g \geq 3$ is not a power of two.
Then there exists a smooth hyperelliptic curve $X$ of genus $g$ so that
$\Jac(X)[p]$ contains $Q$.
\end{corollary}

\begin{proof}
The proof parallels that of Corollary \ref{Cdim2}.
One starts with the supersingular hyperelliptic
curve $X$ of genus 3 and $a$-number $1$ (and thus $\Jac(X)[p] \simeq Q$) from \cite{O}
and works inductively using Proposition \ref{TaddL}.
\end{proof}

It is natural to ask whether Corollary \ref{Cdim2} could be
strengthened to state that $\Jac(X)[p] \simeq N \oplus (\ZZ/p
\oplus \mu_p)^{g-2}$. This raises the following geometric
question.

\begin{question}\label{Qaddmore}
Given any choice of $\Lambda=\{\lambda_1,\ldots,\lambda_{2r}\}$, does there
exist $\mu \in \Aa_k -\Lambda$ so that the hyperelliptic curve branched at
$\{\lambda_1,\ldots,\lambda_{2r},\mu,\infty \}$ is ordinary?
\end{question}

For a generic choice of $\Lambda$, the answer to Question \ref{Qaddmore} is yes by Lemma \ref{Ldegn}.
This question will have an affirmative answer if the hypersurface $D_r$ does not contain any coordinate line
$L(\vec{\lambda}_{2r})$.  The question is equivalent to asking whether,
given a hyperelliptic cover $\phi:X \to \PP_k$, it is always possible
to deform $X$ to an ordinary curve by moving only one branch point.

An affirmative answer to Question \ref{Qaddmore} would allow one to strengthen
Proposition \ref{TaddL} to state that $\Jac(X)[p] \simeq {\mathbb G} \oplus (\ZZ/p \oplus \mu_p)^{g'-1+s}$.
This is because the curves $C_1$ and $C_2$ in the proof can be generically chosen to be
ordinary.
So an affirmative answer to Question \ref{Qaddmore} would imply
that for all $g \geq 4$ there exists a smooth hyperelliptic curve $X$ with
$\Jac(X)[p] \simeq N \oplus (\ZZ/p \oplus \mu_p)^{g-2}$.
If this is true, then $\Jac(X)[p] \simeq N \oplus (\ZZ/p \oplus \mu_p)^{g-2}$ when
$X$ is the generic geometric point of $\CH_g \cap V_{g,g-2}$.

\bibliographystyle{plain}
\bibliography{darren}

\flushright{
Darren Glass\\
Department of Mathematics\\
Columbia University\\
New York, NY 10027\\
glass@math.columbia.edu }

\flushright{
Rachel J. Pries\\
101 Weber Building\\
Colorado State University\\
Fort Collins, CO 80523-1874\\
pries@math.colostate.edu }

\end{document}